\documentclass[12pt]{article}
\usepackage{fullpage, amsmath, amsthm, amssymb}

\newtheorem{theorem}{Theorem}
\newtheorem{lemma}[theorem]{Lemma}
\newtheorem{prop}[theorem]{Proposition}

\theoremstyle{definition}
\newtheorem{defn}[theorem]{Definition}
\newtheorem{example}[theorem]{Example}

\numberwithin{equation}{theorem}

\def\calC{\mathcal{C}}
\def\ZZ{\mathbb{Z}}

\def\dual{\vee}
\def\be{\mathbf{e}}
\def\bv{\mathbf{v}}
\def\bw{\mathbf{w}}

\DeclareMathOperator{\alg}{alg}
\DeclareMathOperator{\Gal}{Gal}
\DeclareMathOperator{\GL}{GL}
\DeclareMathOperator{\perf}{perf}

\title{A formulation of difference Galois theory}
\author{Kiran S. Kedlaya}
\date{October 18, 2006}

\begin{document}

\maketitle

The Galois theory of linear difference equations is somewhat more
complicated than its analogues, the ordinary Galois theory of polynomials
and the Picard-Vessiot theory of linear differential equations.
In a development such as in \cite{singer-vanderput}, 
one restricts to difference fields with algebraically
closed fields of constants. In this note, we point out how a
highly general framework for differential and difference Galois theory
introduced by Andr\'e \cite{andre-diff} makes it possible to generalize
this situation somewhat, at the price of passing from difference fields
to some difference rings. (The fact that such a price is unavoidable is
demonstrated explicitly in \cite[Chapter~1]{singer-vanderput}.)
This is meant in part as an advertisement for
\cite{andre-diff}; thanks to Michael Singer for directing us to
this lovely paper.

\begin{defn}
By a \emph{difference ring/field}, we will mean a commutative
 ring $R$ equipped with an automorphism
$\tau$. (In some texts,
this definition would characterize an \emph{inversive} difference ring/field.)
A \emph{difference module} over a difference ring $R$ is an $R$-module $M$
equipped with a bijective semilinear $\tau$-action (i.e., for
$r \in R$ and $m \in M$, $\tau(rm) = \tau(r) \tau(m)$).
A difference submodule of $R$ itself is called a \emph{difference ideal}.
The \emph{ring of constants}
$C_R$ of $R$ is the subring of $R$ fixed by $\tau$.
\end{defn}

\begin{defn}
A nonzero difference module is \emph{simple} if it contains no nonzero
difference submodules; if a difference ring is simple as a difference
module over itself, we say it is a \emph{simple difference ring}.
A simple difference
ring is necessarily reduced, since the nilradical is a difference ideal.
Also, a localization of a simple difference ring at a
multiplicative subset stable under $\tau^{\pm 1}$
is again simple.
\end{defn}

\begin{lemma} \label{L:total ring}
Let $R$ be a simple difference ring, and let $S$ be the total ring
of fractions of $R$ (i.e., the localization of $R$ at the multiplicative set
of non-zero divisors). Then $C_R$ is a field and $C_S = C_R$.
\end{lemma}
\begin{proof}
For any $x \in C_R$ nonzero, $xR$ is a nonzero
difference ideal and so must coincide with $R$. Hence $x$ has an inverse
in $R$, which must also be fixed by $\tau$.

For $x \in C_S$ nonzero, let $I$ be the set of $y \in R$ 
such that $xy \in R$. Then $I$ is a nonzero difference ideal
of $R$, so must coincide with $R$; that is, $x \in C_S \cap R = C_R$.
\end{proof}

\begin{defn}
Let $R$ be a simple 
difference ring, and let $M$ be a difference
module over $R$ whose underlying $R$-module is finite projective.
A \emph{Picard-Vessiot extension} of $S$ for $M$ is a difference
ring $S$ extending $R$ such that:
\begin{enumerate}
\item[(a)] $S$ is faithfully flat as an $R$-module;
\item[(b)] $S$ is a simple difference ring;
\item[(c)] $C_S = C_R$;
\item[(d)] $M \otimes_R S$ is trivial (i.e., admits a basis fixed by $\tau$);
\item[(e)] $S$ is generated as an $R$-algebra by the coefficients
used to express the elements of
$\tau$-fixed bases of $M$ and $M^\dual$ in terms of generators
of $M$ and $M^\dual$.
\end{enumerate}
Since the description of (e) is a bit terse, we write it out
in symbols. Choose generators $\bv_1, \dots, \bv_m$
and $\bw_1, \dots, \bw_m$ for $M$ and $M^\dual$ respectively.
Choose a $\tau$-fixed basis 
$\be_1, \dots, \be_n$ of $M \otimes_R S$,
and let $\be^*_1, \dots, \be^*_n$ be the basis of $M^\dual \otimes_R S$ dual to
$\be_1, \dots, \be_n$. Write $\be_j = \sum_i A_{ij} \bv_i$ and
$\be^*_j = \sum_i B_{ij} \bw_i$. Then (e) asserts that $S$
is generated as an $R$-algebra by the $A_{ij}$ and the $B_{ij}$.
\end{defn}

\begin{defn}
We say a difference ring is \emph{ind-noetherian} if it
is a direct limit of difference subrings whose underlying rings
are noetherian. A typical difference ring not having this property is 
the ring $\ZZ[x_i: i \in \ZZ]$ with the
automorphism $\tau(x_i) = x_{i+1}$.
\end{defn}

\begin{prop} \label{P:Picard-Vessiot}
Let $R$ be a simple ind-noetherian
difference ring with algebraically closed field of constants.
Let $M$ be a difference module over $R$ whose underlying $R$-module
is finite projective.
Then there exists a unique Picard-Vessiot extension for $M$, and this
difference ring is itself ind-noetherian.
\end{prop}
\begin{proof}
Existence of a Picard-Vessiot extension is a consequence of
\cite[Th\'e\-or\`eme~3.4.3.1]{andre-diff} applied to a suitable noetherian
subring $R_0$ of $R$; the resulting ring $S$ has the form $R \otimes_{R_0} S_0$
for some ring $S_0$ of finite type over $R_0$. Since $S_0$ is noetherian
and $R$ is ind-noetherian, $S$ is ind-noetherian.
Uniqueness of the Picard-Vessiot extension (ind-noetherian
or not)
follows from \cite[Corollaire~3.4.2.5]{andre-diff}.
\end{proof}

\begin{defn}
Let $R$ be a simple ind-noetherian
difference ring with algebraically closed field of constants.
Let $M$ be a difference module over $R$ whose underlying $R$-module
is finite projective, and let $S$ be the Picard-Vessiot extension of $R$ for 
$M$, as produced by Proposition~\ref{P:Picard-Vessiot}. 
Then the tannakian category generated by $M$ admits
a fibre functor over $C_R$
sending $N$ to $(N \otimes_R S)^\tau$, the set of $\tau$-fixed elements of 
$N \otimes_R S$. The automorphism group of this fibre
functor is called the \emph{difference Galois group} of $S$ over $R$, or of $M$
over $R$. It is an affine algebraic group over $C_R$,
and its $C_R$-points are canonically identified with the automorphisms
of $S$ as a difference ring over $R$ \cite[Th\'eor\`eme~3.5.1.1]{andre-diff}.
Moreover, if $C_R$ has characteristic zero, then the fixed subring of $S$
under the action of the difference Galois group
is precisely $R$ \cite[Lemme~3.5.2.1]{andre-diff}.
\end{defn}

We originally got interested in this construction
via the concept of Galois descent,
which can be understood as follows.

\begin{prop} \label{P:Galois descent}
Let $R$ be a simple ind-noetherian
difference ring with algebraically closed field of constants.
Let $S$ be a Picard-Vessiot extension of $R$ with difference Galois group $G$.
Let $V$ be a finite dimensional $C_R$-vector space, and let $\rho: G \to 
\GL(V)$ be an algebraic representation. 
Then $V$ occurs as the action of $G$ on $(M \otimes_R S)^\tau$
for some difference module $M$ over $R$ whose base
extension to $S$ is trivial.
\end{prop}
\begin{proof}
This is an immediate consequence of Tannaka duality:
the representation $\rho$ must appear somewhere in
the Tannakian category generated by a difference module $N$ with Picard-Vessiot
extension $S$.
\end{proof}

\begin{defn} \label{D:universal}
Let $R$ be a simple ind-noetherian
difference ring with algebraically closed field of constants.
Let $\calC$ be a tannakian subcategory of the category of all difference
modules over $R$, whose underlying $R$-modules are finite projective.
We then have a functor taking each 
difference module $M$ over $R$ whose underlying $R$-module is 
finite projective to $(M \otimes_R S)^\tau$, where $S$ is a Picard-Vessiot
extension for $M$. This gives a fibre functor over $C_R$; 
define the \emph{difference Galois group} of $\calC$ as the
automorphism group of this functor.
By taking a suitable direct limit, we can construct (noncanonically)
a \emph{universal Picard-Vessiot extension} for $\calC$, upon which the
difference Galois group acts. Again, if $C_R$ has characteristic zero,
the fixed subring of this action will be $R$.
\end{defn}

Finally, we point out what one has gained from working at this level
of generality.
\begin{example}
Recall again that in \cite{singer-vanderput}, one must start with a difference
field with algebraically closed field of constants. Suppose instead that
$R$ is a difference field, or even an ind-noetherian simple difference ring,
with arbitrary field of constants $C_R$. We may then view $S = R \otimes_{C_R}
C_R^{\alg}$ as an ind-noetherian simple difference ring with the action
on $C_R^{\alg}$ being trivial; one easily verifies that
$C_S = C_R^{\alg}$, so $S$ admits a Picard-Vessiot theory as above.
In fact, the action of $\tau$ on $S$ commutes with the action of 
$\Gal(C_R^{\alg}/C_R^{\perf})$, so one can combine
descent as above with ordinary Galois descent on $C_R^{\alg}$,
at least in case $C_R$ is perfect.
\end{example}

\end{document}